\documentclass[twoside,12pt]{article}
\usepackage{amsmath}
\usepackage[mathscr]{eucal}
\usepackage[charter]{mathdesign}

\newtheorem{THEOREM}{Theorem}
\newtheorem{LEMMA}[THEOREM]{Lemma}
\newtheorem{DEFINITION}[THEOREM]{Definition}
\newtheorem{REMARK}[THEOREM]{Remark}
\newtheorem{PROPOSITION}[THEOREM]{Proposition}
\newtheorem{COROLLARY}[THEOREM]{Corollary}

\def\beginabs{\begin{abstract}\noindent}
\def\endabs{\end{abstract}}
\def\beginrem{\begin{REMARK}\rm}
\def\endrem{\end{REMARK}}
\def\beginthm{\begin{THEOREM}}
\def\endthm{\end{THEOREM}}
\def\beginlem{\begin{LEMMA}}
\def\endlem{\end{LEMMA}}
\def\beginprop{\begin{PROPOSITION}}
\def\endprop{\end{PROPOSITION}}
\def\begincor{\begin{COROLLARY}}
\def\endcor{\end{COROLLARY}}
\def\begindef{\begin{DEFINITION}\rm}
\def\enddef{\end{DEFINITION}}

\newcommand{\cab}{,\allowbreak}
\newcommand{\cas}{,\allowbreak\,}
\newcommand{\imply}{\relax\allowbreak\ifmmode\;{\Rightarrow}\allowbreak
\;\else\ifhmode\unskip\enskip\fi$\Rightarrow$\allowbreak\enskip\fi}

\newcommand{\rds}{\abovedisplayskip 1ex plus .3 ex minus .3 ex
\belowdisplayskip 1ex plus .3 ex minus .3 ex}

\newcommand{\bigspace}{\par\ifdim\lastskip<2.5 ex 
\removelastskip\penalty-200 \vskip 2.5 ex plus .3ex minus .3ex\fi}
\newcommand{\fairspace}{\par\ifdim\lastskip<2 ex 
\removelastskip\penalty-150 \vskip 2 ex plus .25ex minus .25ex\fi}
\newcommand{\paraspace}{\par\ifdim\lastskip<1.25 ex 
\removelastskip\penalty-100 \vskip 1.25 ex plus .25ex minus .25ex\fi}
\newcommand{\smallspace}{\par\ifdim\lastskip<1 ex 
\removelastskip\penalty-75 \vskip 1 ex plus .2ex minus .2ex\fi}

\newcommand{\proofspace}{\par\ifdim\lastskip<2 ex
\removelastskip\penalty-150  \vskip 2 ex plus .25ex minus .25ex\fi}

\def\endproofsymbol{\leavevmode
\hbox{\rm\vrule height 1.4ex  width 1.4ex depth .8ex }}
\def\endproofsymbol{\mbox{$\blacksquare$}}

\def\endproof{\Tag{\endproofsymbol}\proofspace}
\newcommand{\proof}{\par\noindent{\bf Proof.}\kern.5em}

\makeatletter
\newcommand{\X}{\relax\ifmmode\mathscr X\else $\m@th\mathscr X$\fi}

\newcommand\lgcover[1]{\relax\ifmmode\left\langle#1\right\rangle
\else$\m@th\left\langle#1\right\rangle$\fi}

\newcommand{\abs}[1]{\relax\ifmmode\left|#1\right|\else
$\m@th\left|#1\right|$\fi}
\newcommand{\cardinal}[1]{\relax\ifmmode{|}#1{|}\else
$\m@th{|}#1{|}$\fi}

\newcommand{\norm}[1]{\relax\ifmmode\left\Vert#1\right\Vert\else
$\m@th\left\Vert#1\right\Vert$\fi}
\def\seq#1#2{\relax\ifmmode#1_1,#1_2\cab\ldots\cab#1_{#2}\else
$\m@th#1_1,#1_2\cab\ldots\cab#1_{#2}$\fi}

\newcommand{\Tag}[1]{\ifvmode\else\unskip\fi
\nobreak\hfil\penalty50 \hskip2em \null
\nobreak\hfil#1\skip@\parfillskip\parfillskip\z@skip
\count@\finalhyphendemerits\finalhyphendemerits\z@\par
\parfillskip\skip@\finalhyphendemerits\count@}

\newdimen\LabeLmargin\LabeLmargin=0pt
\newdimen\LabeLwidth\LabeLwidth=0pt
\newdimen\Hangamount\Hangamount=0pt
\newdimen \ListSpace  \ListSpace=2ex
\newdimen\LabelSpace \LabelSpace=1ex

\newcommand{\setListSpace}{\ifdim\lastskip<\ListSpace\removelastskip
\penalty-100 \vskip\ListSpace plus .15 \ListSpace minus .15 \ListSpace\fi} 

\newcommand{\multilist}[2]{\par\begingroup\setListSpace
\parindent=0pt \LabeLmargin=\Hangamount\setbox0=\hbox{#1}
\advance\LabeLmargin by\wd0 \Hangamount=\LabeLmargin
\setbox0=\hbox{#2}\advance\Hangamount by\wd0 \LabeLwidth=\wd0}

\newcommand{\singlelist}[1]{\par\begingroup\setListSpace
\parindent=0pt \setbox0=\hbox{#1}\LabeLwidth=\wd0 }

\newcommand{\setLabelSpace}{\ifdim\lastskip<\ListSpace
\removelastskip\penalty-75  \vskip \LabelSpace plus .15 \LabelSpace
minus .15 \LabelSpace\fi}

\newcommand{\listitem}[1]{\par\setLabelSpace \hangindent=\LabeLwidth
\hangafter=1 \leavevmode \hbox to\LabeLwidth{#1\hfill}\ignorespaces}

\newcommand{\listend}{\par\setListSpace\endgroup}

\def\@begintheorem#1#2{\trivlist\item[\hskip\labelsep{\bfseries
#1\ #2.}]\itshape} 
\def\@seccntformat#1{\csname the#1\endcsname.\quad}
\makeatother

\begin{document}

\thispagestyle{empty}
\begin{center}
\large\bf A Note on Set Graceful Labeling of Graphs\end{center}

 \begin{center}
{\bf G. R. Vijayakumar}\\
School of Mathematics,
 Tata Institute of Fundamental Research\\
 Homi Bhabha Road, Colaba, Mumbai 400\,005, 
    India\\ Email: vijay@math.tifr.res.in\end{center}


\beginabs We settle affirmatively a conjecture posed
in [S. M. Hegde, Set colorings of graphs,
{\it  European Journal of  Combinatorics\/}  {\bf30} (4)
(2009),  986--995]: If  some subsets of a set $X$ are 
 assigned injectively to all vertices of a complete bipartite
graph $G$ such that the collection of all sets, each of which is
the symmetric difference of the sets assigned to the ends of
some edge, is the set of all nonempty subsets of $X$, then $G$
is a star.
\endabs

\begingroup \ListSpace=1ex \rightskip 0pt plus 1fil
\singlelist{{\bf Keywords:}\enspace}
\listitem{\bf Keywords:}  set graceful labeling of a graph,
bipartite graph, symmetric difference of two sets.
\listend\endgroup


\noindent Infinite graphs and graphs which are not simple are
out of our consideration.
Let $G$ be a graph with vertex set $V$ and edge set $E$ and 
 $\X$ be the power set of a set. Let $f$ be  a function
from  $V$  to $\X$;  we associate with $f$,  a function from
$E$ to $\X$, denoted by  $\hat f$:
 for all $xy\in E$, $\hat f(xy)=f(x)\Delta f(y)$.
 (The symmetric difference of 
two sets $A$ and $B$ is denoted by $A\Delta B$.) 
If both $f$ and $\hat f$ are injective
and the range of the latter is $\X\setminus\{\emptyset\}$,
then $f$ is called a {\it set graceful
labeling\/} of $G$.
 Settling \cite[Conjecture 2]{h} is the objective of this note:
\fairspace \noindent {\bf Theorem.}\enspace
{\sl If a complete bipartite graph $G$ has
a set graceful labeling, then it is a star}. 
\proofspace
\proof Let $V$ and $E$ be respectively, the vertex set and the
edge set of $G$ and $P\cas Q$ be the bipartition of $G$.
Let $f:V\to \X$ be a set graceful labeling of
$G$. Since $\hat f$ is a bijection from $E$
to $\X\setminus \{\emptyset\}$, it
follows that
$\cardinal E=\cardinal{\X}-1$; i.e.,
$$\rds  \cardinal P\cardinal Q+1=\cardinal{\X}.\eqno(1)$$
 Suppose that $G$ is not a star. Then $\cardinal
P\ne 1\ne\cardinal Q$; this implies that $(\cardinal
P-1)(\cardinal Q-1)>0$; i.e., $\cardinal P\cardinal
Q+1>\cardinal P+\cardinal Q$; therefore by (1), $\cardinal \X>\cardinal V$;
since $f: V\to \X$ is injective, we can find some $A\in
\X\setminus \{f(v):v\in V\}$. 
Define a map $g:V\to \X$ as follows: for any
$v\in V$, $g(v)=A\Delta f(v)$. Since for any $uv\in E$,
$g(u)\Delta g(v)=A\Delta f(u)\Delta A\Delta
f(v)=f(u)\Delta f(v)$, it follows that
$\hat g=\hat f$; further $g$ is obviously
injective; therefore, $g$  is a set graceful labeling.
Since $A\notin\{f(v): v\in V\}$, for each $v\in V$, $g(v)=A\Delta
f(v)\ne \emptyset$. Therefore,
$$\rds \emptyset\notin\{g(v):v\in V\}.\eqno (2)$$ \paraspace
Now, let $p\in P$.  Let us show
that there is exactly one vertex $v\in P$
such that $g(p)\Delta g(v)\in \{g(x):x\in Q\}$.
 Since $g(p)\neq \emptyset$ by (2), and
the range of $\hat g$ is $\X\setminus\{\emptyset\}$, 
 for some $p'\in P$ and $q\in Q$,
$g(p)=\hat g(p'q)$; i.e., $g(p)=g(p')\Delta g(q)$;
therefore, $g(p)\Delta g(p')\allowbreak=g(q)$. Suppose that $p''\in
P$ and $r\in Q$ such that $g(p)\Delta
g(p'')=g(r)$. Then $g(p')\Delta g(q)\allowbreak =g(p)=g(p'')\Delta g(r)$;
i.e., $\hat g(p'q)=\hat g(p''r)$. Since
$\hat g$  is injective, 
 $p'=p''$. Thus for any element $p\in P$, there is a unique element
$v\in P$ such that $g(p)\Delta g(v)\in
\{g(x):x\in Q\}$. Define a map $\theta:P\to P$ as follows. For any $v\in
P$, let $\theta(v)$ be the (unique) vertex in $P$ such that $g(v)\Delta
g(\theta(v))\in \{g(x):x\in Q\}$.
Note that by (2), for each $v\in P$,
$\theta(v)\ne v$ and  $\theta^2(v)=v$. Thus $\bigl\{\{v\cab
\theta(v)\}:v\in P\bigr\}$ is a partition of $P$
into subsets of order 2. Therefore $\cardinal P$ is even
whence by (1), $\cardinal\X$ 
is odd---a contradiction. \endproof

\end{document}